\documentclass[12pt]{amsart}
\usepackage{amssymb,latexsym}

\newtheorem{thm}[equation]{Theorem}

\newtheorem{prop}[equation]{Proposition}
\newtheorem{cor}[equation]{Corollary}
\newtheorem{thma}{Theorem}[section]
\newtheorem*{thm*}{Theorem}
\newtheorem*{lem*}{Lemma}
\newtheorem*{prop*}{Proposition}
\newtheorem*{cor*}{Corollary}
\theoremstyle{remark}
\newtheorem{rmk}[equation]{Remark}
\newtheorem*{rmk*}{Remark}

\newcommand{\gl}{\mathfrak{gl}}

\newcommand{\so}{\mathfrak{so}}
\renewcommand{\sp}{\mathfrak{sp}}
\newcommand{\GL}{\mathsf{GL}}

\newcommand{\End}{\operatorname{End}}

\newcommand{\g}{\mathfrak{g}}
\newcommand{\h}{\mathfrak{h}}
\newcommand{\N}{{\mathbb N}}
\newcommand{\Z}{{\mathbb Z}}
\newcommand{\Q}{{\mathbb Q}}

\newcommand{\E}{{\mathsf E}}
\newcommand{\Schur}[1]{\mathsf{S}{(#1)}}
\newcommand{\U}{{\mathfrak U}}
\newcommand{\epsi}{\varepsilon}
\newcommand{\al}{\alpha}

\newcommand{\divided}[2]{#1^{(#2)}}
\renewcommand{\ge}{\geqslant}
\renewcommand{\le}{\leqslant}
\numberwithin{equation}{subsection}
\allowdisplaybreaks

\begin{document}
\title[{\sf Presenting generalized Schur algebras}]%
{Presenting generalized Schur Algebras in Types $B$, $C$, $D$}

% author 1 information
\author[\sf Doty]{Stephen Doty}
\address{Department of Mathematics and Statistics\\ 
Loyola University Chicago\\ 
Chicago, Illinois 60626 USA}
%\curraddr{}
\thanks{The authors are grateful to Ms.\ Qunhua Liu of Tsinghua 
University for pointing out an error in a preliminary version 
of this manuscript.}
\email{doty@math.luc.edu}

% author 2 information
\author[\sf Giaquinto]{Anthony Giaquinto} 
\address{Department of Mathematics and Statistics\\ 
Loyola University Chicago\\ 
Chicago, Illinois 60626 USA}
%\curraddr{}
\email{tonyg@math.luc.edu}

% author 3 information
\author[\sf Sullivan]{John Sullivan} 
\address{Department of Mathematics, University of Washington, Seattle,
Washington 98195-4350}
\thanks{The third author thanks the Mathematics Departments 
of Bowdoin College and Yale University for their generous support and 
hospitality.}
\email{sullivan@math.washington.edu}

\date{January 2005}

\begin{abstract}
We give explicit presentations by generators and relations of certain
generalized Schur algebras (associated with tensor powers of the
natural representation) in types $B$, $C$, $D$.  This extends previous
results in type $A$ obtained by two of the authors. The presentation
is compatible with the Serre presentation of the corresponding
universal enveloping algebra. In types $C$, $D$ this gives a
presentation of the corresponding classical Schur algebra (the image
of the representation on a tensor power) since the classical Schur
algebra coincides with the generalized Schur algebra in those types.
This coincidence between the generalized and classical Schur algebra
fails in type $B$, in general.
\end{abstract}

\maketitle

\section*{Introduction}

Throughout the paper we work over $\Q$.  Vector spaces are over $\Q$
unless we say otherwise. In fact, all the results are valid over any
field of characteristic zero.

In \cite{PGSA} a new presentation of a generalized Schur algebra was
given.  This presentation is compatible with Lusztig's modified form
$\dot{\U}$ of the universal enveloping algebra $\U$, in which the
Cartan generators are replaced by a system of orthogonal idempotents
corresponding to weight space projectors.  We are interested in the
generalized Schur algebras $S(\pi)$ associated to the set $\pi=\pi(r)$
of dominant weights in a given tensor power $\E^{\otimes r}$ of the
natural representation $\E$ of a simple Lie algebra of type $B$, $C$,
or $D$. Our main result is a presentation by generators and relations
of $S(\pi)$ which is directly compatible with Serre's presentation of
$\U$; i.e., the generators of the zero part are Cartan generators
instead of idempotents. These results are formulated in Section
\ref{sec:main} and proved in Section \ref{sec:pf}; they extend results
of \cite{DG:rank1, DG:PSA} for type $A$ to the other classical types.

We are also interested in the classical Schur algebras $\Schur{r}$,
which we define to be the image of the representation $\U \to
\End(\E^{\otimes r})$. These algebras are closely linked to classical
invariant theory, and interest in them goes back to Schur and Weyl.
They are in types $B$--$D$ the commuting algebras for the action on
tensors of an appropriate Brauer algebra; see \cite{Brauer}. (In type
$A$ the classical Schur algebras are the commuting algebras for the
natural action of symmetric groups, acting by place permutation.)  As
explained in \cite[\S7.4]{PGSA}, the classical Schur algebra coincides
with the corresponding generalized Schur algebra $S(\pi)$ if and only
if the set $\pi_0$ of highest weights labeling composition factors of
$\E^{\otimes r}$ coincides with the set $\pi$ of dominant weights
occurring in $\E^{\otimes r}$. The inclusion $\pi_0 \subset \pi$ is
obvious. In types $A$, $C$, and $D$ it turns out that this inclusion
is equality; in type $B$ this is not so --- see \ref{prop:GSA} and
\ref{rmk:GSA}. Thus the Schur algebra $\Schur{r}$ coincides with the
generalized Schur algebra $S(\pi)$ in types $A$, $C$, and $D$, but in
type $B$ the Schur algebra $\Schur{r}$ is in general a proper quotient
of $S(\pi)$.

We rely on Weyl for the computation of the set $\pi_0$ of highest
weights labeling composition factors of tensor powers of $\E$.  In
order to make this exposition somewhat self-contained, we collect the
relevant results from Weyl's book in an appendix.

Although we work in characteristic zero throughout the paper, it
should be noted that the Schur algebras and generalized Schur algebras
considered herein are defined over $\Z$, so their study can be
undertaken in any characteristic.

Another interesting problem is to find a basis of $S(\pi)$ which is in
some sense compatible with a Poincare-Birkhoff-Witt basis of $\U$.  As
an application of Littelmann's path model, we describe such a basis in 
Section \ref{sec:basis}.

\section{Generalized Schur algebras} \label{sec:GSA}

Let $\g$ be a reductive Lie algebra. The theory of generalized Schur
algebras was introduced by Donkin in \cite{SA1, SA2, SA3}. A
generalized Schur algebra is a certain quotient of $\U=\U(\g)$
obtained by throwing away all but finitely many simple modules. More
precisely, let $\pi$ be a saturated set of dominant weights, meaning
that if $\lambda \in \pi$ and if $\mu$ is a dominant weight such that
$\mu \le \lambda$ (in the usual dominance order), then $\mu \in \pi$.
The generalized Schur algebra determined by $\pi$ is the algebra
$S(\pi):= \U(\g)/\mathcal{I}$, where $\mathcal{I}$ is the ideal of
$\U(\g)$ consisting of all elements of $\U(\g)$ annihilating every
simple $\U$-module of highest weight belonging to $\pi$.

In this paper, we take $\g$ to be a simple Lie algebra of classical
type $B_n$, $C_n$, or $D_n$ and we always take $\pi$ equal to the set
$\Pi^+(\E^{\otimes r})$ of dominant weights of $\E^{\otimes r}$, where
$\E$ is the natural representation of $\g$.  Donkin \cite{SA2,
Donkin:Arcata} showed that, in types $A$ and $C$, the generalized
Schur algebra $S(\pi)$ determined by this choice of $\pi$ coincides
with the Schur algebra $\Schur{r}$. We extend this to type $D$; see
Proposition \ref{prop:GSA} ahead. The corresponding result is not
generally true in type $B$.

\subsection{Basic notation}

In types $B_n,C_n,D_n$ let $\g$ be defined by the form given by
$$
\begin{pmatrix} 0&I&0\\I&0&0\\0&0&1\end{pmatrix};\quad
\begin{pmatrix} 0&I\\-I&0\end{pmatrix};\quad
\begin{pmatrix} 0&I\\I&0\end{pmatrix}.
$$
where $I$ is the $n \times n$ identity matrix.

Let $m = \dim \E$. Set $e_{i,j} = (\delta_{ik}\delta_{lj})_{1\le k,l
\le m}$. The set $\{ e_{i,j} \mid 1\le i,j \le m \}$ is a basis of
$\gl_m$.  In type $A_{n-1}$ we have $m=n$. The $H_i := e_{i,i}$ ($1
\le i \le n$) form a basis for the diagonal Cartan subalgebra $\h$ of
$\gl_n$ In types $B_n$, $C_n$, $D_n$ we have respectively $m = 2n+1,
2n, 2n$ and we set $H_i = e_{i,i} - e_{n+i,n+i}$ ($1 \le i \le n$);
the $H_i$ form a basis for the diagonal Cartan subalgebra $\h$ of
$\so_{2n+1}$, $\sp_{2n}$, $\so_{2n}$ respectively.

In all types $A$--$D$, let $\{ \epsi_j \}$ be the basis of $\h^*$ dual
to the basis $\{ H_i \}$ of $\h$; so that $\epsi_j(H_i) =
\delta_{ij}$. Define a bilinear form $(\ ,\ )$ on $\h^*$ such that
$(\epsi_i, \epsi_j) = \delta_{ij}$.

Denote by $\alpha_1, \dots, \alpha_{n}$ a fixed choice of simple roots
in types $B_n, C_n, D_n$.  Let $(a_{ij})$ be the Cartan matrix,
defined by $a_{ij} = 2(\alpha_i, \alpha_j)/(\alpha_i, \alpha_i)$.

\subsection{Weights}

We regard weights as $n$-tuples of rational numbers, determining
linear functionals on $\h$ by recording their values on the basis
$H_1, \dots, H_n$. In other words, we identify the linear combination
$\lambda_1 \epsi_1 + \cdots + \lambda_n \epsi_n$ with the tuple $\lambda
= (\lambda_1, \dots, \lambda_n)$.

Fix simple root vectors $e_i \in \g_{\alpha_i}$, $f_i \in
\g_{-\alpha_i}$.  The fundamental weights $\varpi_j$ ($1 \le j \le n$)
are the elements of $\h^*$ defined by $\varpi_j(h_i) = \delta_{ij}$,
where $h_i:= [e_i,f_i]$. The lattice of integral weights is the free
abelian group $X$ generated by the $\varpi_i$, and the set $X^+$ of
dominant weights is the cone $\sum \N\varpi_i$.  One checks that $X$
is generated by $\epsi_1, \dots, \epsi_n$ in type $C_n$ and by
$\epsi_1, \dots, \epsi_n$ together with the element $(\epsi_1 + \cdots
+ \epsi_n)/2$ in types $B_n$, $D_n$. As usual, the dominance (partial)
order on $X$ is defined by declaring that $\lambda \le \mu$ (for
$\lambda, \mu \in X$) if $\mu-\lambda \in \N\alpha_1 + \cdots +
\N\alpha_n$.

We identify $\h^*$ with $\Q^n$ by regarding $\epsi_1,\dots, \epsi_n$ as the
standard basis of $\Q^n$; this identifies $X$ with a subgroup of
$\Q^n$. In fact, under this identification, we have
\begin{equation}
X = \begin{cases}
\Z^n & \text{(type $C_n$)};\\
\Z^n \cup ((\frac{1}{2},\dots,\frac{1}{2})+\Z^n) &
\text{(types $B_n$, $D_n$)}.
\end{cases}
\end{equation}
The fundamental weights in type $C_n$ are given explicitly by the
equalities
\begin{equation}
\varpi_i = \epsi_1 + \cdots + \epsi_i \quad (1\le i \le n).
\end{equation}
In type $B_n$ the fundamental weights are given by
\begin{equation}
\varpi_i = \epsi_1 + \cdots + \epsi_i \quad (1\le i \le n-1),
\qquad \varpi_{n} = (\epsi_1 + \cdots + \epsi_{n})/2
\end{equation}
and in type $D_n$ by
\begin{equation}
\begin{gathered}
\varpi_i = \epsi_1 + \cdots + \epsi_i \quad (1\le i \le n-2),\\
\varpi_{n-1} = (\epsi_1 + \cdots +\epsi_{n-1} - \epsi_{n})/2, \quad
\varpi_{n} = (\epsi_1 + \cdots + \epsi_{n})/2.
\end{gathered}
\end{equation}
The set $X^+$ of dominant weights is the set of all $\lambda =
(\lambda_1, \dots, \lambda_n) \in X$ satisfying
\begin{equation}\label{domdef}
\begin{aligned}
\lambda_1 \ge \lambda_2 \ge \cdots \ge \lambda_n \ge 0 \qquad
& \text{(types $B_n$, $C_n$)};\\
\lambda_1 \ge \lambda_2 \ge \cdots \ge \lambda_{n-1} \ge |\lambda_n| \quad
& \text{(type $D_n$)}.
\end{aligned}
\end{equation}

\subsection{Signed compositions}

We shall need the following notations. Write
$$
\begin{aligned}
\Lambda(n,r) &= \{(\lambda_1, \dots, \lambda_n) \in \N^n \mid
\textstyle\sum \lambda_i = r \};\\
\overline{\Lambda}(n,r) &= \{(\lambda_1, \dots, \lambda_n) \in \Z^n \mid
\textstyle\sum |\lambda_i| = r \}.
\end{aligned}
$$
The first set is the set of $n$-part compositions of $r$ and the second is
the set of $n$-part {\em signed} compositions of $r$. Note that we allow 0
to appear in a composition.  Set
$$
\Lambda^+(n,r) = \{(\lambda_1, \dots, \lambda_n) \in \N^n \mid
\textstyle\sum \lambda_i = r,\ \lambda_1 \ge \lambda_2 \ge \cdots \ge
\lambda_n \}.
$$
As usual, we identify members of $\Lambda^+(n,r)$ with partitions
of not more than $n$ parts.  Given $\lambda \in \Lambda^+(n,r)$, set
$\lambda^- := (\lambda_1, \dots, \lambda_{n-1}, -\lambda_n)$ (its
associated weight), and let
$$
\Lambda^-(n,r) = \{\lambda^- \mid \lambda \in \Lambda^+(n,r) \},
\qquad \Lambda^\pm(n,r) = \Lambda^-(n,r) \cup \Lambda^+(n,r).
$$
Label the finite-dimensional simple $\g$-modules $L$ by their highest
weight $\lambda \in X^+$ (regarded as a vector of $\h$-eigenvalues
on $H_1, \dots, H_n$).

\begin{prop} \label{p1}
(a) The set of weights $\Pi$ of $\E^{\otimes r}$ is the set of all
signed $n$-part compositions of $r-2j$ for $0\le j \le [r/2]$ (i.e.\
the union of the $\overline{\Lambda}(n,r-2j)$ for $0\le j \le [r/2]$)
in types $C_n$, $D_n$ and the set of all signed $n$-part compositions
of $r-j$ for $0\le j \le r$ (i.e.\ the union of the
$\overline{\Lambda}(n,r-j)$ for $0\le j \le r$) in type $B_n$.

(b) The set $\Pi^+$ of dominant weights of $\E^{\otimes r}$ is the
union over $0\le j \le [r/2]$ of the sets $\Lambda^+(n,r-2j)$ in type
$C_n$, the union over $0\le j \le r$ of the sets $\Lambda^+(n,r-j)$ in
type $B_n$, and the union over $0\le j \le [r/2]$ of the sets
$\Lambda^\pm(n,r-2j)$ in type $D_n$.
\end{prop}

\begin{proof}
In types $C_n, D_n$ the weights of $\E$ are $\{ \pm \epsi_1, \dots,
\pm \epsi_n \}$. In type $B_n$ the weights of $\E$ are $\{ \pm
\epsi_1, \dots, \pm \epsi_n \} \cup \{0\}$. The weights
of $\E^{\otimes r}$ are given by all expressions of the form
$$
w_{1} + \cdots +  w_{r}
$$
where $w_i$ is a weight of $\E$ for each $i$. (The $w_i$ are not
necessarily distinct.) Part (a) is now clear.

To prove part (b), combine part (a) with \eqref{domdef}.
\end{proof}

\begin{prop} \label{p2}
The set $\pi = \Pi^+(\E^{\otimes r})$ of dominant weights in
$\E^{\otimes r}$ is a saturated subset of $X^+$, for types $B_n$,
$C_n$, $D_n$.
\end{prop}

\begin{proof}
One can decompose $\E^{\otimes r}$ into a direct sum of irreducible
  modules. The set of dominant weights of an irreducible is
  necessarily saturated (in characteristic zero), and the union of
  saturated sets is necessarily saturated. In, fact, this argument
  shows that the set of dominant weights of any $\g$-module must be a
  saturated set.
\end{proof}

We remark that one can also give a combinatorial proof of the
preceding result.

\begin{prop} \label{prop:GSA}
The Schur algebras $\Schur{r}$ in types $C_n$, $D_n$ are generalized
Schur algebras determined by the saturated set $\pi =
\Pi^+(\E^{\otimes r})$.
\end{prop}

\begin{proof}
$\Schur{r}$ is by definition the image of the representation $\U \to
\End(\E^{\otimes r})$, so $\Schur{r} \simeq \U/A$ where $A$ is the
annihilator of $\E^{\otimes r}$.  By Wedderburn theory, the simple
$\Schur{r}$-modules are the direct summands of $\E^{\otimes r}$, and
thus must have highest weight belonging to the set $\pi =
\Pi^+(\E^{\otimes r})$.

Let $S(\pi) = \U/\mathcal{I}$ be Donkin's generalized Schur algebra
determined by the saturated set $\pi$. Here $\mathcal{I}$ is the
ideal of $\U$ consisting of the elements annihilating every simple
$\U$-module of highest weight belonging to $\pi$. Clearly $A \subseteq
\mathcal{I}$.

Let $\pi_0$ be the set of highest weights of composition factors
appearing as a direct summand in a Wedderburn decomposition of
$\E^{\otimes r}$.  Weyl computed the decomposition of tensor space for
the special orthogonal and symplectic groups. His results show that,
in types $C_n$ and $D_n$, $\pi_0 = \pi$. (See the appendix for a
detailed summary of Weyl's results, with references.) This justifiess
the equlity $A = \mathcal{I}$, in types $C$, $D$. The proof is
complete.
\end{proof}

\begin{rmk} \label{rmk:GSA}
The preceding result often fails for type $B_n$.  Indeed, the natural
module $\E$ has two dominant weights but only one highest weight, so
for any $n$ the needed equlity $\pi_0 = \pi$ fails already for $r=1$.
Moreover, in type $B_2 = \so_5$ we have $\pi_0 = \{(2,0), (1,1),
(0,0)\}$ but the set $\pi$ is $\{(2,0), (1,1), (1.0), (0,0)\}$. Here
we relied on Weyl's Theorem \ref{dog} in the appendix for a
description of $\pi_0$ in type $B$, and Proposition \ref{p1} above for
the set $\pi$.

Interestingly, the desired equality $\pi_0=\pi$ holds when $n=1$ and
$r \ge 2$.  It would be useful to classifiy those pairs $n$ and $r$ for
which $\pi_0=\pi$ for type $B_n$. \hfill$\diamond$
\end{rmk}

\subsection{The idempotent presentation} \label{IdempotentPresentation}

{\setcounter{equation}{0}
\renewcommand{\theequation}{R\arabic{equation}} Let $\pi =
\Pi^+(\E^{\otimes r})$, in types $B$, $C$, and $D$.  In
\cite[6.13]{PGSA}, it was shown that the generalized Schur algebra
$S(\pi)$ is isomorphic with the associative algebra (with 1) on
generators $e_i$, $f_i$ ($1\le i \le n$), $1_\lambda$ ($\lambda \in
W\pi$) with the relations
\begin{gather}
1_\lambda 1_\mu = \delta_{\lambda\mu} 1_\lambda, \quad
\sum_{\lambda\in W\pi} 1_\lambda = 1 \\
e_i f_j - f_j e_i = \delta_{ij} \sum_{\lambda\in W\pi}
(\alpha_i^\vee, \lambda)\, 1_\lambda \\
e_i 1_\lambda =
\begin{cases}
1_{\lambda+\alpha_i} e_i &
   \text{if $\lambda+\alpha_i \in W\pi$}\\
0 & \text{otherwise}
\end{cases} \\
f_i 1_\lambda =
\begin{cases}
1_{\lambda-\alpha_i} f_i &
   \text{if $\lambda-\alpha_i \in W\pi$}\\
0 & \text{otherwise}
\end{cases} \\
1_\lambda e_i =
\begin{cases}
e_i 1_{\lambda-\alpha_i} &
   \text{if $\lambda-\alpha_i \in W\pi$}\\
0 & \text{otherwise}
\end{cases} \\
1_\lambda f_i =
\begin{cases}
f_i 1_{\lambda+\alpha_i} &
   \text{if $\lambda+\alpha_i \in W\pi$}\\
0 & \text{otherwise}
\end{cases} \\
\sum_{s=0}^{1-a_{ij}} (-1)^s \binom{1-a_{ij}}{s}
e_i^{1-a_{ij}-s}e_je_i^s = 0 \quad (i \ne j) \\
\sum_{s=0}^{1-a_{ij}} (-1)^s \binom{1-a_{ij}}{s}
f_i^{1-a_{ij}-s}f_jf_i^s = 0 \quad (i \ne j).
\end{gather}
}

\par\noindent
Here $W$ is the Weyl group attached to the Lie algebra $\g$, the
$a_{ij}$ are as before, and $\alpha_i^\vee = 2\alpha_i/(\alpha_i,
\alpha_i)$ for $i=1, \dots, n$. Note that $W\pi$ is equal to
$\Pi(\E^{\otimes r})$, the set described explicitly in Proposition 
\ref{p1}.

\section{Main results}\label{sec:main}

In the statements to follow, notice that the first, third, fourth,
fifth, and sixth relations are identical in all types.  In other
words, only the second and seventh relations vary by type. (The
seventh relation is the same in types $C$, $D$.)

\subsection{Type $B$}

The root system for $B_n$ is realized by $\alpha_i = \epsi_i -
\epsi_{i+1}$ for $i<n$; $\alpha_n = \epsi_n$.

\begin{thm}\setcounter{equation}{0}\label{t1}
\renewcommand{\theequation}{B\arabic{equation}} 
Over $\Q$, the generalized Schur algebra $S(\pi)$ of type $B_n$ is
isomorphic with the associative algebra (with 1) on generators $e_i$,
$f_i$, $H_i$ ($1\le i \le n$) and with relations
\begin{gather}
H_iH_j = H_jH_i \\
e_if_j - f_je_i = \begin{cases}
   \delta_{ij}(H_i-H_{i+1}) & (i<n) \\
   \delta_{ij}(2H_n) & (i=n)
   \end{cases} \\
H_ie_j - e_jH_i = (\varepsilon_i, \alpha_j)e_j, \quad
H_if_j - f_jH_i = -(\varepsilon_i, \alpha_j)f_j\\
\sum_{s=0}^{1-a_{ij}} (-1)^s \binom{1-a_{ij}}{s}
e_i^{1-a_{ij}-s}e_je_i^s = 0 \quad (i \ne j) \\
\sum_{s=0}^{1-a_{ij}} (-1)^s \binom{1-a_{ij}}{s}
f_i^{1-a_{ij}-s}f_jf_i^s = 0 \quad (i \ne j) \\
(H_i + r)(H_i+r-1)(H_i+r-2)\cdots(H_i - r)=0\\
(J+r)(J+r-1)(J+r-2)\cdots(J-r+1)(J-r)=0
\end{gather}
where $J = \pm H_1 \pm H_2 \pm \cdots \pm H_n$ varies over all $2^n$
possible sign choices.
\setcounter{equation}{1}
\end{thm}

Note that the enveloping algebra $\U(\so_{2n+1})$ is the algebra on the
same generators but subject only to the relations (B1)--(B5);
moreover, that presentation of $\U(\so_{2n+1})$ is equivalent to the
usual Serre presentation.

The relation (B6) is necessary.

\subsection{Type $C$}

The root system for $C_n$ is realized by $\alpha_i = \epsi_i -
\epsi_{i+1}$ for $i<n$; $\alpha_n = 2\epsi_n$.

\begin{thm}\setcounter{equation}{0}\label{t2}
\renewcommand{\theequation}{C\arabic{equation}} Over $\Q$, the
generalized Schur algebra $S(\pi)$ of type $C_n$ (which coincides with
the Schur algebra $\Schur{r}$) is isomorphic with the associative
algebra (with 1) on generators $e_i$, $f_i$, $H_i$ ($1\le i \le n$)
and with relations
\begin{gather}
H_iH_j = H_jH_i \\
e_if_j - f_je_i = \begin{cases}
   \delta_{ij}(H_i-H_{i+1}) & (i<n) \\
   \delta_{ij}H_n & (i=n)
   \end{cases} \\
H_ie_j - e_jH_i = (\varepsilon_i, \alpha_j)e_j, \quad
H_if_j - f_jH_i = -(\varepsilon_i, \alpha_j)f_j\\
\sum_{s=0}^{1-a_{ij}} (-1)^s \binom{1-a_{ij}}{s}
e_i^{1-a_{ij}-s}e_je_i^s = 0 \quad (i \ne j) \\
\sum_{s=0}^{1-a_{ij}} (-1)^s \binom{1-a_{ij}}{s}
f_i^{1-a_{ij}-s}f_jf_i^s = 0 \quad (i \ne j) \\
(H_i + r)(H_i+r-1)(H_i+r-2)\cdots(H_i - r)=0\\
(J + r)(J+r-2)(J+r-4)\cdots(J - r+2)(J - r)=0
\end{gather}
where $J = \pm H_1 \pm H_2 \pm \cdots \pm H_n$ varies over all
possible sign choices. 
\setcounter{equation}{2}
\end{thm}

Note that the enveloping algebra $\U(\sp_{2n})$ is the algebra on the
same generators but subject only to the relations (C1)--(C5); that
presentation of $\U(\sp_{2n})$ is equivalent to the usual Serre
presentation.

The relation (C6) is superfluous.

\subsection{Type $D$}

The root system for $D_n$ is realized by taking $\alpha_i = \epsi_i -
\epsi_{i+1}$ for $i<n$; $\alpha_n = \epsi_{n-1}+\epsi_n$.

\begin{thm}\setcounter{equation}{0}\label{t3}
\renewcommand{\theequation}{D\arabic{equation}} Over $\Q$, the
generalized Schur algebra $S(\pi)$ of type $D_n$ (which coincides with
the Schur algebra $\Schur{r}$) is isomorphic with the associative
algebra (with 1) on generators $e_i$, $f_i$, $H_i$ ($1\le i \le n$)
and with relations
\begin{gather}
H_iH_j = H_jH_i \\
e_if_j - f_je_i = \begin{cases}
   \delta_{ij}(H_i-H_{i+1}) & (i<n) \\
   \delta_{ij}(H_{n-1}+H_n) & (i=n)
   \end{cases} \\
H_ie_j - e_jH_i = (\varepsilon_i, \alpha_j)e_j, \quad
H_if_j - f_jH_i = -(\varepsilon_i, \alpha_j)f_j\\
\sum_{s=0}^{1-a_{ij}} (-1)^s \binom{1-a_{ij}}{s}
e_i^{1-a_{ij}-s}e_je_i^s = 0 \quad (i \ne j) \\
\sum_{s=0}^{1-a_{ij}} (-1)^s \binom{1-a_{ij}}{s}
f_i^{1-a_{ij}-s}f_jf_i^s = 0 \quad (i \ne j) \\
(H_i + r)(H_i+r-1)(H_i+r-2)\cdots(H_i - r)=0\\
(J + r)(J+r-2)(J+r-4)\cdots(J - r+2)(J - r)=0
\end{gather}
where $J = \pm H_1 \pm H_2 \pm \cdots \pm H_n$ varies over all
possible sign choices.
\end{thm}

The enveloping algebra $\U(\so_{2n})$ is the algebra on the
same generators but subject only to the relations (D1)--(D5); that
presentation of $\U(\so_{2n})$ is equivalent to the usual Serre
presentation.

The relation (D6) is superfluous.

\bigskip

The proof of all three theorems of this section is given in the next
section. Our strategy is to show that the presentation of the theorem
is equivalent to the idempotent presentation of section
\ref{IdempotentPresentation}.

\begin{rmk}
One can easily show that any $H_i$, viewed as an operator on
$\E^{\otimes r}$, satisfies its minimal polynomial $P_1(T)$, and
similarly that any $J$, viewed as an operator on $\E^{\otimes r}$,
satisfies its minimal polynomial $P_1(T)$, in type $B$, or $P_2(T)$,
in types $C$, $D$.
\hfill$\diamond$\end{rmk}

\section{Proof of the main theorems}\label{sec:pf}

We will show that the generalized Schur algebra $S(\pi)$, with $\pi =
\Pi^+(\E^{\otimes r})$, defined by the presentation in
\ref{IdempotentPresentation} is isomorphic with the algebra given by
the generators and relations of Theorems \ref{t1}, \ref{t2}, or
\ref{t3}, in types $B_n$--$D_n$.

\subsection{The algebra $\Phi$}

Write $\U = \U(\g)$.  Given a positive integer $r$, set
\begin{equation}\label{P1P2}
\begin{aligned}
P_1(T) &= (T+r)(T+r-1)\cdots (T-r+1)(T-r),\\
P_2(T) &= (T+r)(T+r-2)\cdots (T-r+2)(T-r)
\end{aligned}
\end{equation}
polynomials of degree $2r+1$, $r+1$, respectively.  Let $\Phi$ be the
algebra given by the generators and relations of Theorem \ref{t1},
\ref{t2}, or \ref{t3}.  Then $\Phi = \U/I$. In types $C_n$ and $D_n$,
$I$ is the two-sided ideal of $\U$ generated by the $P_1(H_i)$ ($i=1,
\dots, n$) and $P_2(J)$ for all $J=\pm H_1 \pm \cdots \pm H_n$. In
type $B_n$, $I$ is the two-sided ideal of $\U$ generated by the
$P_1(H_i)$ ($i=1, \dots, n$) and $P_1(J)$ for all $J$.

From the triangular decomposition $\U = \U^- \U^0 \U^+$ of $\U$ we
have a corresponding triangular decomposition $\Phi = \Phi^- \Phi^0
\Phi^+$, where each algebra $\Phi^-$, $\Phi^0$, $\Phi^+$ is defined to
be the image under the appropriate surjective map of the corresponding
subalgebra of $\U$.  Let $\U_\Z$ be Kostant's $\Z$-form of $\U$
relative to the Chevalley generators $e_i$, $f_i$; this is the
$\Z$-subalgebra of $\U$ generated by all $\divided{f_i}{a}$,
$\divided{e_i}{c}$ ($a, c \in \N$, $1 \le i \le n$).  Then we have
equalities $\U_\Z = \U_\Z^- \U_\Z^0 \U_\Z^+$, $\Phi_\Z = \Phi_\Z^-
\Phi_\Z^0 \Phi_\Z^+$ where the various subalgebras are defined
in the obvious manner.

For the moment, regard $H_1, \dots, H_n$ as commuting
indeterminates. A given set of polynomials in the polynomial ring
$\Q[H_1, \dots, H_n]$ determines an affine variety in $\Q^n$.

\begin{prop}
Let $V \subseteq \Q^n$ be the common zero locus of $P_2(J)$ for all $J
= \pm H_1 \pm H_2 \pm \cdots \pm H_n$, in types $C_n$, $D_n$.  In type
$B_n$ let $V \subseteq \Q^n$ be the common zero locus of $P_2(J)$ for
all $J$, along with $P_1(H_i)$ for all $i=1, \dots, n$.  Then $V =
\Pi$, the set of weights of $\E^{\otimes r}$.
\end{prop}

\begin{proof}
For fixed $r$, let $T = \{-r,-r+2,\ldots ,r-2,r\}$. This is the set of
$m \in \Z$ satisfying $|m| \le r$, $m \equiv r \pmod{2}$.  The set
$\Pi$, in types $C_n$ and $D_n$, is the set of $(\lambda_1, \ldots,
\lambda_n) \in \Z^n$ satisfying $\sum |\lambda_i| \in T$.

Consider $v=(v_1,\ldots ,v_n)\in V$, a solution to $P_2(J)=0$ for all
choices of $J=\pm H_1 \pm \cdots \pm H_n$. Then $\pm v_1\pm
v_2\pm\cdots \pm v_n\in T$ for all possible sign choices. In
particular, for every $i$ with $1\le i\le n$ we have
$$
v_1+\cdots+v_n\in T \text{ and } -v_1-\cdots
-v_{i-1}+v_i-v_{i+1}-\cdots -v_n\in T.
$$
From this we conclude that $2v_i$ is an even integer since the sum
of any two elements of $T$ is even. Thus $v_i\in \Z$ for all $i$.
Since $\sum |v_i|\in T$, the inclusion $V\subseteq \Pi$ holds.

Now take $\lambda =(\lambda_1,\ldots ,\lambda_n)\in \Pi$. We need to
prove that $\pm \lambda_1+\cdots \pm \lambda_n\in T$ for all sign
choices. By the description of $\Pi$ given above, we know $\sum
|\lambda_i|\in T$. For any choice $\sigma_i$ of signs, we have the
congruence $\sum \sigma_i|\lambda_i| \equiv \sum |\lambda_i|
\pmod{2}$. Moreover, $\left|\sum \sigma_i |\lambda_i|\right|\le \sum
|\lambda_i| \le r$.  Thus $\Pi \subseteq V$. This proves that $V=\Pi$
in types $C_n$ and $D_n$.

Now we turn to type $B_n$. Let $T'=\{-r,-r+1,\ldots r-1,r\}$ and note
that in this case $\Pi$ is the set of $\lambda_1,\ldots ,\lambda_n)\in
\Z^n$ satisfying $\sum |\lambda_i| \in T'$. By an argument similar to
the above, it follows that the set of solutions to the equations
$P_1(J)=0$ for all $J$ coincides with the set of $(v_1,\ldots, v_n)\in
(\frac{1}{2}\Z)^n$ such that $\sum |v_i|\in T'$. If the additional
equations $P_1(H_i)=0$ are imposed, then it is clear that the solution
set is reduced exactly to $\Pi = \Pi(\E^{\otimes r})$. The proof is
complete.
\end{proof}

\begin{rmk}
The proof shows, in particular, that relations (C6), (D6) are
consequences of relations (C7), (D7).
\hfill$\diamond$\end{rmk}

We now consider the algebra $\Phi^0=\U^0/(\U^0\cap I)$. By the PBW
theorem, $\U^0$ is isomorphic with the algebra $\Q[H_1, \dots, H_n]$
of polynomials in commuting ``indeterminates'' $H_1, \dots,
H_n$. Define an algebra $\Phi'=\U^0/I^0$ where $I^0$ is the ideal in
$\U^0$ generated by $P_2(J)$ for all $J$ and $P_1(H_i)$ for all $i$,
in types $C_n$, $D_n$, and is the ideal in $\U^0$ generated by
$P_1(J)$ for all $J$ and by $P_1(H_i)$ for all $i$, in type $B_n$.
(In types $C_n$, $D_n$ the generators $P_1(H_i)$ are not needed.)
Given $\lambda \in \Pi$, define an element $1_\lambda \in \Phi^0$ by
\begin{equation}\label{idempdef}
1_\lambda = \prod_i
\frac{P^{(\lambda_i)}_1(H_i)}{P^{(\lambda_i)}_1(\lambda_i)}
\end{equation}
where $P_1^{(k)}(T)$ equals $P_1(T)$ with factor $(T-k)$ deleted, for
a given $k$ satisfying $-r \le k \le r$.  Since
$(H_i-\lambda_i)P_1^{(\lambda_i)}(H_i) = 0$ (by definition of $\Phi$)
we see from \eqref{idempdef} that 
\begin{equation}\label{H_itimes}
H_i 1_\lambda = \lambda_i 1_\lambda \qquad (\lambda \in \Pi, 1\le i
\le n).
\end{equation}

\begin{prop}
(a) The algebra $\Phi^0$ is isomorphic with the algebra $\Phi'$.

(b) The set of all $1_\lambda$ ($\lambda \in \Pi$) is a $\Q$-basis for
$\Phi^0$ and a $\Z$-basis for $\Phi^0_\Z$; moreover, this set is a set
of pairwise orthogonal idempotents in $\Phi^0$ which add up to $1$.
\end{prop}

\begin{proof}
View the $H_i$ as coordinate functions on $\Q^n$. The algebra
$\Phi'$ is the ring of regular functions on the variety of common
zeros of $I^0$. By the previous proposition, this variety is the
finite set $\Pi$. The coordinate ring of $\Pi$ is just the product
$\prod_{\lambda \in \Pi} \Q_{\lambda}$ where $\Q_{\lambda}\cong \Q$
is the function ring of the $\lambda.$ (The only functions defined
at a single point are the constants.) The explicit
isomorphism $\Phi^0 \cong \prod_{\lambda \in \Pi} \Q_{\lambda}$
is realized by the map (denoted by $\phi$) which sends $f(H_1,
\ldots, H_n)$ to $\prod_{\lambda \in \Pi} f(\lambda_1,\ldots, \lambda_n)$.
It is easy to check that $1_{\lambda}(\lambda') =
\delta_{\lambda\lambda'}$. Thus, $\phi (1_\lambda)$
is a vector whose entries are all zero except for one entry which equals one.
Since $\phi$ is an isomorphism, it follows that the set $1_\lambda$
($\lambda \in \Pi$) is a $\Q$-basis for $\Phi'$ and this set is a
set of pairwise orthogonal idempotents of $\Phi'$ which add up to
$1$.

By the definition of $\Phi$ we have an algebra surjection $\U \to
\Phi$ By restriction, this induces an algebra surjection $\U^0 \to
\Phi^0$. The canonical quotient map $\U \to \U/I$ induces, upon
restriction to $\U^0$, a map $\U^0 \to \U/I=\Phi$. The image of this
map is $\Phi^0$ and its kernel is $\U^0 \cap I$, so $\Phi^0 \cong
\U^0/(\U^0 \cap I)$. Clearly $I^0 \subseteq \U^0 \cap I$.  Thus we
obtain an algebra surjection:
\begin{equation} \label{seq}
\Phi' = \U^0/I^0 \rightarrow \Phi^0.
\end{equation}
The dimension of $\Phi'$ is the cardinality of $\Pi$, the set of
weights appearing in the representation $E^{\otimes r}$. 

We consider the quotient of the polynomial ring $\Q[H_1, \dots, H_n]$
(the variables $H_i$ commute) by the ideal $I^0$. It suffices to show
that this quotient is isomorphic with $|\Pi|$ copies of the base field
$\Q$.  By the Chinese remainder theorem, applied repeatedly to the
factors of the polynomial $P_1(H_i)$, for each $i = 1, \dots, n$, we
obtain an isomorphism
$$
\Q[H_1, \dots, H_n]/I^0 \simeq \prod_{\lambda_1, \dots, \lambda_n} 
\Q[H_1, \dots, H_n]/(H_1-\lambda_1, \dots, H_n-\lambda_n, P_k(J)) 
$$ 
where $k = 1$ in types $C_n$, $D_n$ and $k=2$ in type $B_n$. In the
product, each integer $\lambda_i$ belongs to the interval
$[-r,r]$. Each factor in the product is either $\Q$ or zero because
the relations $H_1-\lambda_1 = 0, \dots, H_n - \lambda_n = 0$ make
each variable $H_i$ a constant. Consider the factor for a selection of
constants $\lambda_1, \dots, \lambda_n$. If those values satisfy the
identity $P_k(J) = 0$ for every choice of $J$, then the selection of
constants gives a weight of the $r$th tensor power of $\E$, and the
factor $\Q[H_1, \dots, H_n]/(H_1-\lambda_1, \dots, H_n-\lambda_n,
P_k(J))$ is $\Q$. If the selected values do not satisfy the relation
$P_k(J) = 0$ for every choice of $J$, then the evaluation of such a
relation at $\lambda_1, \dots, \lambda_n$ gives a nonzero constant in
the ideal, so the quotient $\Q[H_1, \dots, H_n]/(H_1-\lambda_1, \dots,
H_n-\lambda_n, P_k(J))$ is 0.

It follows that $\Phi^0$ has the same dimension as $\Phi'$. Thus the
surjection \eqref{seq} is an isomorphism of algebras.  This proves
assertions (a) and (b).
\end{proof}

\begin{rmk}\label{12}
(a) The argument shows in particular that $\U^0 \cap I = I^0$, an
equality which is not obvious from the definitions.

(b) From the proposition and \eqref{H_itimes} it follows immediately (by
multiplication by $1 = \sum_\lambda 1_\lambda$) that in $\Phi$ we have
the equality $H_i = \sum_{\lambda\in \Pi} \lambda_i 1_\lambda$ for any
$i$.
\hfill$\diamond$\end{rmk}

\begin{prop}
The elements $e_i$, $f_i$ ($1\le i \le n$), $1_\lambda$ ($\lambda \in
\Pi$) of $\Phi$ satisfy relations (R1)--(R8), with $\Pi = W\pi$.
\end{prop}

\begin{proof}
Relation (R1) was proved in the previous proposition. To prove (R2),
first consider the case $i \ne j$. Then relations (B2), (C2), (D2) all
assert that $e_if_j-f_je_i = 0$, which is precisely relation (R2) in
this case. Now suppose $i=j$ is strictly less than $n$. Then
$\alpha_i^\vee = \alpha_i = \epsi_i-\epsi_{i+1}$ for types $B_n$,
$C_n$, and $D_n$. Thus $(\alpha_i^\vee,\lambda) = \lambda_i -
\lambda_{i+1}$. Since $i<n$, relations (B2), (C2), (D2) all assert
that $e_i f_i - f_i e_i=H_i-H_{i+1}$. By Remark \ref{12}(b),
$$\textstyle
H_i-H_{i+1}=\sum_{\lambda\in \Pi} (\lambda_i-\lambda_{i+1})1_{\lambda}.
$$
This proves (R2) in the case $i=j < n$.  Now consider the final
remaining case $i=j=n$. In type $B_n$, we have $\alpha_n^{\vee}
=2\alpha_n=2\epsi_n$, and $e_nf_n-f_ne_n=2H_n$. Relation (R2) for type
$B_n$ now follows since $2H_n1_{\lambda}=(\alpha_n^{\vee},
\lambda)1_{\lambda}$ for all $\lambda.$. In type $C_n$,
$\alpha_n^{\vee} =\alpha_n/2=\epsi_n$ and $e_n f_n - e_n f_n=H_n$.
Relation (R2) for type $C_n$ now follows since
$H_n1_{\lambda}=(\alpha_n^{\vee}, \lambda)1_{\lambda}$ for all
$\lambda.$. Finally, in type $D_n$, $\alpha_n^\vee =\alpha_n
=\epsi_{n-1}+\epsi_n$ and $e_n f_n - e_n f_n=H_{n-1}+H_n$. In exactly
the same way as for the other cases, it follows at once that (R2)
holds for type $D_n$.

We now prove relation (R3). First, from relations (B3), (C3), and (D3)
we see that $e_jH_i=(H_i-(\epsi_i,\al_j))e_j$ and so from
\eqref{idempdef} we obtain the equality (in $\Phi$)
\begin{equation} 
e_j 1_{\lambda} = P(H_1,\ldots , H_n)\,e_j
\end{equation}
where $P(H_1,\ldots , H_n)$ is defined by
\begin{equation} 
P(H_1,\ldots , H_n) = \prod_i
\frac{P_1^{(\lambda_i)}(H_i-(\epsi_i,\al_j))}
{P_1^{(\lambda_i)}(\lambda_i)}.
\end{equation}
From Remark \ref{12}(b) and the definition of $P(H_1, \ldots, H_n)$ we
obtain the equality
\begin{equation}
P(H_1,\ldots , H_n)= \textstyle\sum_{\mu\in \Pi}P(\mu_1, \ldots,
\mu_n)1_{\mu}.
\end{equation}
In order to analyze this expression, first note that by its definition
$P_1^{(\lambda_i)}(x)=0$ for all integers $x$ except $x=\lambda_i$ or
$|x|\ge r+1$. Thus for a given $\mu\in \Pi$, $P(\mu_1,\ldots ,\mu_n)
= 0$ unless, for all $i\in \{1,\ldots,n\}$, one of the conditions
$\mu_i-(\epsi_i,\al_j)=\lambda_i$ or $|\mu_i-(\epsi_i,\al_j)|\ge r+1$
holds.

Now suppose there exists $\mu\in \Pi$ and $i\in \{1,\ldots, n\}$ such
that $|\mu_i-(\epsi_i,\al_j)|\ge r+1.$ In this case, even though
$P(\mu_1,\ldots, \mu_n)$ need not vanish, the product $1_{\mu}e_j$ is
necessarily zero as we now show.  We have the equality
$1_{\mu}(H_i-\mu_i)=0$, so trivially $1_{\mu}(H_i-\mu_i)e_j=0$. But,
using relation (B3), (C3), or (D3), we may rewrite the last equality
in the form
\begin{equation}\label{invertible}
1_{\mu}e_j(H_i+(\epsi_i,\al_j)-\mu_i)=0.
\end{equation}
The rightmost factor, $H_i+(\epsi_i,\al_j)-\mu_i$ can be expressed via
Remark \ref{12}(b) as the sum
$$
H_i+(\epsi_i,\al_j)-\mu_i = \sum_{\nu \in \Pi} (\nu_i
+(\epsi_i,\al_j)-\mu_i) 1_\nu
$$ 
in which every coefficient differs from zero since we are under the
assumption that $|\mu_i-(\epsi_i,\al_j)| \ge r+1$ and each
component of an element of $\Pi$ lies in the interval $[-r,r]$. It
follows that $H_i+(\epsi_i,\al_j)-\mu_i$ is an invertible element of
$\Phi$ and so we can multiply equation \eqref{invertible} by its
inverse on the right to obtain the desired result that $1_{\mu}e_j=0$.

Thus, the only $\mu\in \Pi$ for which $P(\mu_1,\ldots
\mu_n)1_{\mu}e_j\neq 0$ is determined by
$\mu_i-(\epsi_i,\al_j)=\lambda_i$ for all $i\in\{1,\ldots , n\}$. Thus
$\mu=\lambda+\al_j$.  Moreover, one easily sees that $P(\lambda+\al_j)
= 1$. Relation (R3) now follows.  The proofs of relations (R4)--(R6)
are similar. Finally, (R7) and (R8) hold since these are among the
defining relations for $\U(\g)$. The proof is complete.
\end{proof}

We note the following corollary for later reference.

\begin{cor}\label{14}
With $\pi = \Pi^+(\E^{\otimes r})$, there is a surjective algebra
homomorphism $S(\pi) \to \Phi$ given by $e_i \to e_i$, $f_i \to f_i$,
$1_\lambda \to 1_\lambda$.
\end{cor}

\subsection{The algebra $S$}

Let $S = S(\pi)$ be the generalized Schur algebra, given by the
generators and relations (R1)--(R8) of \ref{IdempotentPresentation},
for types $B_n$--$D_n$.  Let $\Pi = \Pi(\E^{\otimes r}) = W\pi$, the
set of weights of $\E^{\otimes r}$.  Define elements $H_i \in S$ by
\begin{equation}\label{Hidef}
H_i = \sum_{\lambda \in \Pi} (\epsi_i, \lambda) 1_\lambda =
\sum_{\lambda \in \Pi} \lambda_i 1_\lambda.
\end{equation}

\begin{prop}
With $H_i$ as above, the elements $H_i$, $e_i$, $f_i$ in $S$ satisfy
the relations (B1)--(B7), (C1)--(C7), (D1)--(D7) in types $B_n$, $C_n$,
$D_n$ respectively.
\end{prop}

\begin{proof}
The first, third, fourth, fifth, and sixth relations are the same for
all types, so the argument differs only for the second and seventh
relations. Moreover, the fourth and fifth relations are the same as
(R7) and (R8), so we only need to establish the first, second, third,
sixth, and seventh relations.

It follows from the definition \eqref{Hidef} of the elements $H_i$ and
the commutativity of the $1_\lambda$ that the elements $H_i$, $H_j$
commute in $S$.  This proves the first relation (B1), (C1), (D1).

We consider the third relation. We will show that $H_ie_j - e_jH_i =
(\epsi_i, \alpha_j)e_j$. At this point, it is convenient to set
$1_\lambda = 0$ for all $\lambda\in X - \Pi$. Then the sums in (R1),
(R2), and \eqref{Hidef} can be taken over $X$.  The relations (R3),
(R5) may be expressed by the single equality
\begin{equation} \label{R'3}
e_i 1_\lambda = 1_{\lambda +\alpha_i} e_i  \qquad
(\text{all } i,\  \lambda \in X)
\end{equation}
and (R4), (R6) may be expressed as
\begin{equation} \label{R'4}
f_i 1_\lambda = 1_{\lambda-\alpha_i} f_i \qquad
(\text{all } i,\  \lambda \in X).
\end{equation}
From \eqref{Hidef}, \eqref{R'3}, by reindexing the first sum we obtain
$$
\begin{aligned}
H_ie_j - e_jH_i &= \textstyle\sum_\lambda (\epsi_i, \lambda) 1_\lambda
e_j - \sum_\lambda (\epsi_i,\lambda )e_j 1_\lambda \\
 &= \textstyle \sum_\lambda (\epsi_i,\lambda )e_j1_{\lambda -\alpha_j}-
\sum_\lambda (\epsi_i,\lambda )e_j1_\lambda \\
 &= \textstyle\sum_\lambda (\epsi_i,\lambda +\alpha_j)e_j 1_\lambda -
\sum_\lambda (\epsi_i,\lambda )e_j1_\lambda \\
 &= \textstyle \sum_\lambda (\epsi_i,\alpha_j)e_j1_\lambda =
(\epsi_i,\alpha_j)e_j,
\end{aligned}
$$
by the second part of (R1), where the sums are taken over $X$. This
proves the first part of (B3), (C3), and (D3). The proof of the second
part is entirely similar, using \eqref{R'4} instead of \eqref{R'3}.

We consider the sixth relation. We have equalities
$$
\begin{aligned}
(H_i&+r)\cdots H_i\cdots (H_i-r) \\
&= (\sum_{\lambda \in \Pi} (\epsi_i,\lambda) 1_\lambda +r) \cdots
(\sum_{\lambda \in \Pi} (\epsi_i,\lambda) 1_\lambda) \cdots
(\sum_{\lambda \in \Pi }(\epsi_i,\lambda) 1_\lambda - r) \\
&= (\sum_{\lambda \in \Pi} ((\epsi_i,\lambda)+r) 1_\lambda) \cdots
(\sum_{\lambda \in \Pi}(\epsi_i,\lambda )1_\lambda) \cdots
(\sum_{\lambda \in \Pi }((\epsi_i,\lambda )-r) 1_\lambda) \\
&= \sum_{\lambda \in \Pi }[((\epsi_i,\lambda )+r)\cdots
(\epsi_i,\lambda )\cdots ((\epsi_i,\lambda )-r)]1_\lambda
\end{aligned}
$$
(since the $1_\lambda$ are orthogonal idempotents).  The last
expression above equals 0, i.e., all its coefficients equal 0, because
for $\lambda \in \Pi$, $(\epsi_i,\lambda) = \lambda_i$ is an integer
between $-r$ and $r$. This proves (B6), (C6), and (D6).

We consider (B2). We must show that $e_if_i-f_ie_i =
\delta_{ij}(H_i-H_{i+1})$, when $i<n$ and $e_if_i-f_ie_i =
\delta_{ij}(2H_n)$, when $i=n$. When $i\ne j$, $e_if_i-f_ie_i = 0$ by
(R2), so we are reduced to the case $i=j$.
By \eqref{Hidef}, for $i<n$ we have
$$
\begin{aligned}
H_i-H_{i+1}
&= \textstyle\sum_{\lambda \in \Pi} (\epsi_i-\epsi_{i+1},\lambda )1_\lambda\\
&= \textstyle\sum_{\lambda \in \Pi }(\alpha_i,\lambda )1_\lambda
= \textstyle\sum_{\lambda \in \Pi }(\alpha_i^\vee ,\lambda )1_\lambda
= e_if_i-f_ie_i,
\end{aligned}
$$
by relation (R2). (We used the equality $\alpha_i = \alpha_i^\vee$
for $i<n$.)  For $i=n$, we have
$2H_n=\sum_{\lambda \in \Pi }(2\epsi_n,\lambda
)1_\lambda $ = $\sum_{\lambda \in \Pi }(\alpha_n^\vee ,\lambda
)1_\lambda $ = $e_nf_n-f_ne_n$, by (R2). This proves (B2). The proof
of (C2) and (D2) is similar.

We consider (B7). Consider $J = \sum_{i=1}^n \sigma_i H_i$
for a given choice of signs $(\sigma_1, \dots, \sigma_n)$ in $\{1,
-1\}^n$. Then we have
$$
J = \textstyle\sum_{i=1}^n \sigma_i \sum_{\lambda \in \Pi}
(\epsi_i,\lambda)1_\lambda = \sum_{i=1}^n \sum_{\lambda \in \Pi}
\sigma_i \lambda_i 1_\lambda .
$$
For any integer $s$ we have equalities
$$
\begin{aligned}
(J+s)&= (\textstyle \sum_i \sum_\lambda \sigma_i \lambda_i 1_\lambda)
+ (\sum_\lambda s 1_\lambda) \\
&= \textstyle \sum_\lambda ((\sum_i \sigma_i \lambda_i) + s)1_\lambda
\end{aligned}
$$
where $i$ varies from 1 to $n$ and $\lambda$ varies over $\Pi$.  Hence
we obtain
$$
\begin{aligned}
(J+r)\cdots J \cdots (J-r)
&= \prod_{s=-r}^r \sum_{\lambda\in \Pi}
((\sum_{i=1}^n \sigma_i \lambda_i) + s)1_\lambda\\
&= \sum_{\lambda\in \Pi} \prod_{s=-r}^r
((\sum_{i=1}^n \sigma_i \lambda_i) + s) 1_\lambda.
\end{aligned}
$$
The last expression vanishes since, for each $\lambda \in \Pi $, $\pm
\lambda_1 \pm \cdots \pm \lambda_n$ is an integer between $-r$ and
$r$.  This proves (B7). The proof of (C7), (D7) is similar.
\end{proof}

\begin{cor}\label{16}
There is a surjective algebra homomorphism $\Phi \to S$ mapping
$e_i \to e_i$, $f_i \to f_i$, and $H_i \to H_i= \sum_{\lambda
\in \Pi} (\epsi_i, \lambda) 1_\lambda$.
\end{cor}

\subsection{Conclusion of the proof}

Corollary \ref{14} shows that $S$ is a quotient of $\Phi$. Corollary
\ref{16} shows that $\Phi$ is a quotient of $S$. It follows that
$\Phi$ is isomorphic with $S = S(\pi)$, and hence Theorems \ref{t1},
\ref{t2}, and \ref{t3} are proved.

\newcommand{\ff}{\tilde{f}}
\newcommand{\ee}{\tilde{e}}
\section{A basis for $S(\pi)$}\label{sec:basis}

Let $\g$ be a Lie algebra of classical type defined over $\Q$.  As
before, let $\pi = \Pi^+(\E^{\otimes r})$. We would like to give a
canonical basis for the generalized Schur algebra $S(\pi)$ in terms of
the elements of the universal enveloping algebra.

Since $S(\pi)$ is the direct sum of endomorphism algebras of simple
factors of $\E^{\otimes r}$, we consider first the problem of finding
a basis for the full matrix algebra $\End (M)$ of an irreducible
module $M$. To give a basis for $\End(M)$ in terms of the elements of
the universal enveloping algebra, we will exploit Littelmann's basis
for $M$.

\subsection{Littelmann's basis}

Let $\h$ be a Cartan subalgebra of $\g$. Let $R\subset \h^*$ be the set
of roots of $\h$ in $\g$, i.e., the nonzero eigenvalues for the adjoint
representation of $\h$ on $\g$, and let $S$ be the set of simple roots
relative to some hyperplane in $\h^*$. For each simple root $\alpha $,
let $s_\alpha $ be the corresponding element of the Weyl group of the
root system.  Fix a reduced expression $s_{\alpha_1}\cdots
s_{\alpha_m}$ for the longest word $w_0$ in the Weyl group of $\g$ in
terms of the set $\{s_\alpha \mid \alpha\in S\}$.

Let $X_{\Q }$ be the rational span of the weight lattice of $\g$
within $\h^*$. Consider the set of paths $x:[0,1]\to X_{\Q }$ that
begin at the origin, and take $\Pi $ to be the free $\Z $-module on
that set.  In \cite[\S1]{L1}, Littelmann defines certain
operators $\{\ff_\alpha \}_{\alpha\in S}$ and $\{\ee_\alpha
\}_{\alpha\in S}$ on $\Pi $. The value of $\ff_\alpha $ at a path
$x(t)$ with endpoint $x(1)$ either is the $0$-element of $\Pi $ or is
a particular path $\ff_\alpha \cdot x$ with endpoint $x(1)-\alpha
$. In an inverse sense, the value of $\ee_\alpha $ at a path $x(t)$
with endpoint $x(1)$ either is the $0$-element of $\Pi $ or it is a
particular path $\ee_\alpha \cdot x$ with endpoint $x(1)+\alpha $.
Let $\mathbf{A}$ be the algebra generated by those operators.

Consider any path $P_\lambda $ in the dominant chamber that terminates
at an integral weight $\lambda $.  $P_\lambda $ generates an
irreducible $\mathbf{A}$-module $M_\lambda $. In terms of the chosen
reduced expression $w_0=s_{\alpha_1}\cdots s_{\alpha_m}$, consider the
elements of $\mathbf{A}$ of the form $\ff_{\alpha_1}^{n_1}\cdots
\ff_{\alpha_m}^{n_m}$. Littelmann shows that the elements
$\ff_{\alpha_1}^{n_1}\cdots \ff_{\alpha_m}^{n_m}\cdot P_\lambda $ span
$M_\lambda $. Moreover, in \cite[\S\S6--7]{L2}, he gives a geometric
description of a set $S_\lambda $ of sequences of exponents
$n_1,...,n_m$ such that, as $(n_1,...,n_m)$ ranges over $S_\lambda $,
the elements $\ff_{\alpha_1}^{n_1}\cdots \ff_{\alpha_m}^{n_m}\cdot
P_\lambda $ form a basis for $M_\lambda $.

For each simple root $\alpha $, let $f_\alpha $ and $e_\alpha $ be
nonzero elements of the root spaces of $\g$ corresponding respectively
to the roots $-\alpha $ and $\alpha$, satisfying $[e_\alpha, f_\alpha]
= \alpha^\vee \in \h^*$. Let $L(\lambda)$ be the irreducible
$\g$-module of highest weight $\lambda$, and $v_\lambda$ be a highest
weight vector.  In \cite[Theorem 10.1]{L2}, using the theory of
crystal bases, Littelmann shows that as the sequences $(n_1,...,n_m)$
range over $S_\lambda $, the elements $f_{\alpha_1}^{n_1}\cdots
f_{\alpha_m}^{n_m}\cdot v_\lambda $ form a basis for $L(\lambda) $.

Let $S_\lambda^{opp}$ be the set of sequences $(n_m,...,n_1)$, where
$(n_1,...,n_m)$ ranges over $S_\lambda $. Let $u_\lambda $ be a lowest
weight vector in $L(\lambda) $.  By the same theory of crystal bases,
as $(t_1,...,t_m)$ ranges over $S_\lambda^{opp}$, the elements
$e_{\alpha_1}^{t_1}\cdots e_{\alpha_m}^{t_m}\cdot u_\lambda $, form a
basis for $L(\lambda) $.

\subsection{Basis for ${\End}(L(\lambda) )$}

Consider the dual module $L(\lambda)^*$. Its weights are the negatives
of the weights of $L(\lambda) $, and its highest weight is
$-w_0(\lambda )$, since the lowest weight of $L(\lambda) $ is
$w_0(\lambda )$. Furthermore, the lowest weight of $L(\lambda)^*$ is
$-\lambda $, and as a lowest weight vector, we can take the element
$\chi_{v_\lambda }$ that is 1 at $v_\lambda $ and 0 at all weight
vectors of other weights. Hence, by the preceding paragraph, as
$(t_1,...,t_m)$ ranges over $S_{-w_0(\lambda )}^{opp}$, the elements
$e_{\alpha_1}^{t_1}\cdots e_{\alpha_m}^{t_m}\cdot\chi_{v_\lambda }$,
form a basis for $L(\lambda)^* $.

Let $1_\lambda$ be the element of $\End(L(\lambda))$ that acts as $1$ on
the highest weight line of $L(\lambda)$ and acts as $0$ on the other
weight spaces.

\begin{thm}\label{thm:basis}
Let $\lambda $ be a dominant integral weight. As
$\{n_1,...,n_m\}$ range over $S_\lambda $, and as $\{t_1,...,t_m\}$
range over $S_{-w_0(\lambda )}^{opp}$, the elements
$$
\{f_{\alpha_1}^{n_1}\cdots f_{\alpha_m}^{n_m}\, 1_\lambda \,
e_{\alpha_1}^{t_1}\cdots e_{\alpha_m}^{t_m}\}
$$ 
form a basis for $\End(L(\lambda) )$.
\end{thm}

\begin{proof}  
Under the natural identification of $L(\lambda)\otimes
L(\lambda)^*$ with $\End(L(\lambda) )$, the element 
$$
f_{\alpha_1}^{n_1}\cdots f_{\alpha_m}^{n_m}\cdot v_\lambda\otimes
e_{\alpha_1}^{t_1}\cdots e_{\alpha_m}^{t_m}\cdot \chi_{v_\lambda }
$$
of $L(\lambda)\otimes L(\lambda)^* $ corresponds to the element
$$ 
f_{\alpha_1}^{n_1}\cdots f_{\alpha_m}^{n_m}\, 1_\lambda \,
e_{\alpha_1}^{t_1}\cdots e_{\alpha_m}^{t_m}
$$ 
of $\End(L(\lambda) )$.
\end{proof}

\begin{cor} \label{cor:basis}
In Lie algebras of type $B_n$, $C_n$, or $D_n$, as
$\{n_1,...,n_m\}$ range over $S_\lambda $, and as $\{t_1,...,t_m\}$
range over $S_{\lambda }^{opp}$, the elements
$$
\{f_{\alpha_1}^{n_1}\cdots f_{\alpha_m}^{n_m}\, 1_\lambda \,
e_{\alpha_1}^{t_1}\cdots e_{\alpha_m}^{t_m}\}
$$ 
form a basis for $\End(L(\lambda) )$.
\end{cor}

\begin{proof}
In those types, $w_0=-I$.
\end{proof}

\subsection{Basis for $S(\pi)$}

Consider the generalized Schur algebra $S(\pi)$ for $\pi =
\Pi^+(\E^{\otimes r})$ in types $B_n$, $C_n$, and $D_n$.  $S(\pi)$ is
isomorphic with the direct sum $\bigoplus_{\lambda \in \pi} 
\End( L(\lambda) )$. In \cite[6.10]{PGSA} it is proved that an idempotent
$1_\mu$ ($\mu \in W\pi = \Pi(\E^{\otimes r})$) acts on any
$S(\pi)$-module $M$ as 1 on the eigenspace of $M$ of value $\mu$, and
acts as 0 on the eigenspaces of $M$ of all other values.

\begin{thm}
The algebra $S(\pi)$ has a basis consisting of all elements of the
form
$$
f_{\alpha_1}^{n_1}\cdots f_{\alpha_m}^{n_m}\, 1_{\lambda(i)} \,
e_{\alpha_1}^{t_1}\cdots e_{\alpha_m}^{t_m}
$$
such that 
$$ \{n_1,...,n_m\}\in S_{\lambda} ;\{t_1,...,t_m\}\in S_{\lambda
}^{opp} 
$$
as $\lambda$ varies over $\pi$.
\end{thm}

\begin{proof} (By induction on the partial order on $\pi$.)
Let $\lambda$ be a maximal weight of $\pi$.  The element
$$
f_{\alpha_1}^{n_1}\cdots f_{\alpha_m}^{n_m}\, 1_{\lambda}\,
e_{\alpha_1}^{t_1}\cdots e_{\alpha_m}^{t_m}
$$ 
is zero in $\End( L(\lambda') )$ for each $\lambda'\in \pi$ different
from $\lambda$. Such elements give a basis for the endomorphism
algebra of $L(\lambda)$.  By induction we may assume the basis has
been established as stated for the generalized Schur algebra $S(\pi
\backslash \{\lambda\})$. Note that $\pi \backslash \{\lambda\}$ is
saturated.
\end{proof}

\appendix
\section{Irreducible factors in the $r$th tensor power of the natural
module}

Weyl \cite{Weyl:book} describes the irreducible factors in the $r$th
tensor power of the natural module for a classical group. For
convenience, we summarize here the results from \cite{Weyl:book} 
needed in the paper.

\subsection*{Type $B_n = \mathsf{SO}_{2n+1}$}
Let $T$ be a diagram with row lengths $f_1\ge f_2\ge \cdots f_t>0$,
with $t\le n$. The diagram $T'$ {\it associated to T} has row lengths
$f_1\ge f_2\ge \cdots f_t\ge f_{t+1}\ge\cdots \ge f_{2n+1-t}$, where
$f_{t+1}=f_{t+2}=\cdots =f_{2n+1-t}=1$, i.e., $T'$ is obtained by
adding $2n-2t+1$ boxes to the first column of $T$. The irreducible
$\mathsf{SO}_{2n+1}$-modules associated to $T$ and $T'$ are isomorphic
with highest weight of exponents $f_1\ge f_2\ge \cdots \ge
f_t$. (Weyl, Chapter 7, equations 9.10 and 9.11.) The pairs $\{T,T'\}$
partition Weyl's set of permissible diagrams (Weyl p.\ 155).

The following is Weyl's Theorem 5.7F for the group $\mathsf{SO}_{2n+1}$.

\begin{thma}
For each pair $\{T,T'\}$ of permissible diagrams, the irreducible
module corresponding to $T$ is a factor of the $r$th tensor power of
the natural module iff $\sum_{j=1}^tf_j=m-2k$, for some $k\ge 0$, and
the irreducible module corresponding to $T'$ is a factor of the $r$th
tensor power of the natural module iff $\sum_{j=1}^{2n+1-t}f_j=m-2k$
for some $k\ge 0$.
\end{thma}

\begin{rmk*}
Weyl's theorem is stated for the full orthogonal group, but in section
9 of Chapter 6, he shows that irreducible modules for the full
orthogonal group remain irreducible for the proper orthgonal group
when the dimension of the natural module is odd.
\end{rmk*}

Because $f_j=1$ for $j=t+1,...,2n+1-t$, we can write the condition on
$T'$ in the theorem as $\sum_{j=1}^tf_j+(2n-2t+1)=m-2k$ for some $k\ge
0$, or equally, as $\sum_{j=1}^tf_j+(2n-2t)=m-2k-1$ for some $k\ge 0$.

We can give the irreducible factors of the $r$th tensor power in the
following theorem:

\begin{thma}\label{cat} 
The irreducible factors of the $r$th tensor power of the natural
 module for $\mathsf{SO}_{2n+1}$ are those whose highest characters
 have exponents $f_1\ge f_2\ge \cdots \ge f_t>0$ with $t\le n$, where
 either

(i) $\sum_{j=1}^tf_j=m-2k$, for some $k\ge 0$, or 

(ii) $\sum_{j=1}^tf_j+(2n-2t)=m-2k-1$, for some $k\ge 0$.
\end{thma}

The theorem can be restated in a slightly different form:

\begin{thma}\label{dog}  
The irreducible factors of the $r$th tensor power of the natural
 module for $\mathsf{SO}_{2n+1}$ are those whose highest characters
 have exponents $f_1\ge f_2\ge \cdots \ge f_n\ge 0$, where either

(i) $\sum_{j=1}^nf_j=m-2k$, for some $k\ge 0$, or

(ii) $\sum_{j=1}^nf_j=m-2k'-1$, for some $k'\ge 0$, for which
$f_{n-k'}\ne 0$.
\end{thma}

\begin{proof}
Set $k'=n+k-t$ to interchange parts (ii) of the theorems \ref{cat} and
\ref{dog}. Note that the condition $k\ge 0$ translates into the
condition $t \ge n-k'$.
\end{proof}

\subsection*{Type $C_n = \mathsf{Sp}_{2n}$}
Consider the diagrams $T$ with row lengths $f_1\ge f_2\ge\cdots \ge
f_n\ge 0$.  For the symplectic group $\mathsf{Sp}_{2n}$, the
irreducible modules correspond to the diagrams $T$. The irreducible
$\mathsf{Sp}_{2n}$-module corresponding to the diagram $T$ is a factor
of the $r$th tensor power of the natural module if and only if
$f_1+f_2+\cdots +f_n=m-2k$ for some $k\ge 0$.

In terms of highest weights, we have the following (Weyl,
Chapter 6, p.\ 175 and Theorem 7.8D):

\begin{thma} Consider sequences of integers $f_1\ge f_2\ge\cdots\ge
f_n\ge 0$. The irreducible factors of the $r$th tensor power of the
natural module for $\mathsf{Sp}_{2n}$ are those whose highest
characters have exponents $(f_1,f_2,...,f_n)$, where $f_1+f_2+\cdots
+f_n=m-2k$, for some $k\ge 0$.
\end{thma}

\subsection*{Type $D_n = \mathsf{SO}_{2n}$}
The {\em permissible} diagrams for the full orthogonal group
$\mathsf{O}_{2n}$ are those whose first two columns have combined
length no more than $2n$. The irreducible modules for the full
orthogonal group $\mathsf{O}_{2n}$ correspond to the permissible
diagrams (Weyl Theorem 5.7F).

The set of permissible diagrams can be partitioned into pairs of
associated diagrams as follows.  Let $T$ be a diagram with row lengths
$f_1\ge f_2\ge\cdots \ge f_t> 0$, with $t\le n$. The diagram $T'$ {\it
associated} to $T$ has row lengths $f_1\ge f_2\ge \cdots\ge f_t\ge
f_{t+1}\ge\cdots \ge f_{2n-t}$, where $f_{t+1}=f_{t+2}=\cdots
=f_{2n-t}=1$, i.e., $T'$ is a permissible diagram obtained from the
permissible diagram $T$ by adding $2n-2t$ boxes to the first column of
$T$. By Weyl's theorem 5.7F for $\mathsf{O}_{2n}$, the irreducible
module corresponding to $T$ is a factor of the $r$th tensor power of
the natural module exactly when $\sum_{j=1}^t f_j=m-2k$, for some
$k\ge 0$, and the irreducible module corresponding to $T'$ is a factor
of the $r$th tensor power of the natural module exactly when
$\sum_{j=1}^{2n-t}f_j=m-2k'$, for some $k'\ge 0$. Because
$f_{t+1}=f_{t+2}=\cdots =f_{2n-t}=1$, that condition can be written as
$\sum_{j=1}^tf_j=m-2(n-t)-2k'$, for some $k'\ge 0$.

Consider the pair $\{T,T'\}$. If $T'$ satisfies the condition
$\sum_{j=1}^tf_j=m-2(n-t)-2k'$, for some $k'\ge 0$, then $T$ satisfies
the condition $\sum_{j=1}^tf_j=m-2k$, for the value $k=(n-t)+k'\ge
0$. Hence, if the diagram $T'$ corresponds to a factor of the $r$th
tensor power, so does the diagram $T$.

When we restrict the irreducible $\mathsf{O}_{2n}$-modules
corresponding to $T$ and $T'$ to the proper orthogonal group
$\mathsf{SO}_{2n}$, there are two cases to consider (Weyl, Theorems
5.9A and 7.9). In the first case, $t<n$. There, $T$ and $T'$ are
distinct diagrams that correspond to nonisomorphic
$\mathsf{O}_{2n}$-modules. Upon restriction to $\mathsf{SO}_{2n}$,
those modules become {\it isomorphic} irreducible
$\mathsf{SO}_{2n}$-modules. The highest character of that irreducible
$\mathsf{SO}_{2n}$-module has nonzero exponents $(f_1,f_2,...,f_t)$,
with $t<n$. In the second case, $t=n$. There, $T=T'$, and the
corresponding irreducible $\mathsf{O}_{2n}$-module, upon restriction
to $\mathsf{SO}_{2n}$, splits into two nonisomorphic irreducible
$\mathsf{SO}_{2n}$-modules, one whose highest character has exponents
$(f_1,f_2,...,f_{n-1},f_n)$ and the other whose highest character has
exponents $(f_1,f_2,...,f_{n-1},-f_n)$.

The following theorem gives the isomorphism classes of the irreducible
$\mathsf{SO}_{2n}$-modules in the $r$th tensor power of the natural
module, summing up the conclusions of the preceding paragraphs.

\begin{thma} Consider sequences of integers $f_1\ge f_2\ge\cdots\ge
f_n\ge 0$ such that $\sum_{j=1}^nf_j=m-2k$, for some $k\ge 0$.. The
irreducible modules in the $r$th tensor power of the natural module
for $\mathsf{SO}_{2n}$ are those whose highest characters have
exponents $(f_1, f_2,...,f_{n-1},f_n)$ with $f_n=0$, and those whose
highest characters have exponents $(f_1,f_2,...,f_{n-1},f_n)$ or
$(f_1,f_2,...,f_{n-1},-f_n)$, with $f_n>0$.
\end{thma}

\end{document}